\documentclass{jnmp03c}

\usepackage{amsmath}

\begin{document}
\setcounter{page}{31}

\newcommand{\rmname}[1]
  {\expandafter\newcommand \csname #1\endcsname {{\operatorname{#1}}}}
\rmname{Diff}
\rmname{pr}
\rmname{Vect}

\renewcommand{\evenhead}{P Grozman}
\renewcommand{\oddhead}{Bilinear Invariant Operators on the Symplectic Manifold}

\thispagestyle{empty}

\FistPageHead{1}{\pageref{grozman-firstpage}--\pageref{grozman-lastpage}}{Letter}

\copyrightnote{2001}{P Grozman}

\Name{On Bilinear Invariant Differential Operators Acting on Tensor Fields on
the Symplectic Manifold}\label{grozman-firstpage}

\Author{Pavel GROZMAN}

\Adress{Department of Mathematics, University of Stockholm\\
Roslagsv. 101, Kr\"aftriket hus 6, S-106 91, Stockholm, Sweden\\
E-mail: grozman@matematik.su.se}

\Date{Received May 6, 2000; Revised August 30, 2000;
Accepted October 3, 2000}

\begin{abstract}
\noindent Let $M$ be an $n$-dimensional manifold, $V$ the space of
a representation $\rho: GL(n)\longrightarrow GL(V)$.  Locally, let
$T(V)$ be the space of sections of the tensor bundle with
fiber~$V$ over a sufficiently small open set $U\subset M$, in
other words, $T(V)$ is the space of tensor fields of type $V$ on
$M$ on which the group $\Diff (M)$ of diffeomorphisms of~$M$
naturally acts.  Elsewhere, the author classified the $\Diff
(M)$-invariant differential operators $D: T(V_{1})\otimes
T(V_{2})\longrightarrow T(V_{3})$ for irreducible fibers with
lowest weight. Here the result is generalized to bilinear
operators invariant with respect to the group $\Diff_{\omega}(M)$
of symplectomorphisms of the symplectic manifold $(M, \omega)$. We
classify all first order invariant operators; the list of other
operators is conjectural. Among the new operators we mention a 2nd
order one which determins an ``algebra'' structure on the space of
metrics (symmetric forms) on~$M$.
\end{abstract}

\noindent
Let $\rho $ be a representation of the group $Sp(2m; {\mathbb R} )$
in a $V_{\rho }$.  A {\it tensor field of type} $\rho $ on a
$2m$-dimensional symplectic manifold $M$ is an object $t$ defined in
each local coordinate system $x$, in which the symplectic form is of
the canonical form $\omega = \sum\limits_{i\le m}dx_{i}\wedge
dx_{2m+1-i}$, by the vector $t(x)\in V_{\rho }$, where the collections
of all vectors $t(x)$ are such that the passage to
other coordinates, $y$ (with the same property), is defined by the
formula
\[
t(y(x))=\rho \left (\frac{\partial y(x)}{
\partial x} \right )t(x).
\]
Traditionally (see reviews \cite{K1,K2}) the fibers of the
tensor bundles were considered finite dimensional, but Leites showed
recently \cite{LKW} that on supermanifolds it is natural and fruitful
to consider infinite dimensional fibers: this leads to semi-infinite
cohomology of supermanifolds.  Similar problem for symplectic
manifolds and supermanifolds was not studied yet.

The space of smooth tensor fields of type $\rho $ will be denoted by
$T(\rho )$ or by $T(\lambda )$, where $\lambda =(\lambda _{1}, \dots ,
\lambda _{m})$ is the lowest (or, for finite dimensional
representations, highest, for convenience) weight of the irreducible
representation $\rho$.

In what follows the letters $\rho$, $\sigma$, $\tau$ will denote
irreducible representations of $Sp(2m; {\mathbb R} )$ and letters $\lambda$,
$\mu$, $\nu $ their highest weights.  (We should have concidered lowest
weight only, but in this report we stick to finite
dimensional representations, where it does not matter.)

\medskip

\noindent
{\bf Examples of spaces of tensor fields:}
\begin{itemize}
\topsep0mm
\partopsep0mm
\parsep0mm
\itemsep0mm
\item[a)] $T(0)=C^{\infty }(M)$;
\item[b)] $T(1, 0, \ldots , 0) \cong \Vect\cong \Omega ^{1}$ is the space of
vector fields or (which is the same on any symplectic manifold thanks to
the nondegeneracy of $\omega$) the space of 1-forms on $M$;
\item[c)] $\prod^{r}=T(\underbrace{1, \dots , 1}_{r\text{-many}}, 0, \dots ,
0)$ the space of primitive $r$-forms.
\end{itemize}

\noindent
{\bf Remark 1.} Observe that the spaces of tensor fields traditionally
are understood as $p$ times covariant and $q$ times contravariant ones
($T^p_{q}$) or their subspaces subject to some symmetry conditions.
Such tensors split into the direct sum of irreducible $Sp(2n)$-modules
and the same module can be encountered in distinct $T^p_{q}$ (with
different $p$'s and $q$'s).  For example, the tensors of trivial type,
$T(0)$ can be encountered in $T^0_{0}=C^{\infty }(M)$, as stated
above, and also in $T^2_{0}=C^{\infty }(M)\omega$ and in many other
places.  In this sense the case of symplectomorphisms differs from
the general diffeomorphisms, where each irreducible $GL(n)$-module has
a unique embedding into the tensor algebra.

\medskip

A differential operator $B : T(\rho _1) \otimes T(\rho _2)\otimes
\cdots \otimes T(\rho _n )\longrightarrow T(\tau)$ is called $n$-{\it
ary} ({\it un}ary, {\it bi}nary, etc.  for $n=1, 2$, respectively). 
Such an operator $B$ is called $\Diff(M)$-{\it invariant} if it is
uniquely expressed in all coordinate systems; it is
$\Diff_{\omega}(M)$-{\it invariant} if it is uniquely expressed in all
coordinate systems in which the symplectic form is of the standard
(canonical) form.

\medskip

\noindent
{\bf The unary $\pbf{\Diff_{\omega}(M)}$-invariant differential
operators:}

A~N~Rudakov classified all such operators (\cite{R1,R2}):

{\bf 0-th order}: the multiplication by a scalar;

{\bf 1-st order}: the derivations of the primitive forms $d_{+}:
\prod^{r} \longrightarrow \prod^{r+1}$ and $d_{-}:
\prod^{r+1}\longrightarrow \prod^{r}$ ($0 \le r \le m-1$).  These
operators are compositions of the exterior differential $d : \Omega
^{p}\longrightarrow \Omega ^{p+1}$ and the projection onto the space
of primitive forms; recall that $ \Omega ^{p} = \prod^{p} \oplus
\prod^{p-2}\omega \oplus \prod^{p-4}\omega^2 \oplus\dots$ for $p \le m$ and $\Omega ^{p} \cong \Omega
^{2m-p}$;

{\bf 2-nd order}: $d_{2} = d_{+}\circ
d_{-} : \prod^{r}\longrightarrow \prod^{r}$ ($1 \le r \le m$).

\medskip

\noindent
{\bf Remark 2.}  Rudakov's theorem implies that other invariant
operators that might spring to mind ($d_-\circ \omega \circ d_-$, etc.) are
multiples of the described ones.

\medskip

\noindent
{\bf The binary $\pbf{\Diff_{\omega}(M)}$-invariant differential
operators:}

On the space $T_{c}(\rho )$ of tensor fields of type $\rho$ with
compact support, as indicated by the subscript, there is an invariant
inner product
\renewcommand{\theequation}{{\rm IP}}
\begin{equation}
\langle \chi , \theta \rangle = \int_{M}
\langle \chi (x), \theta (x)\rangle \omega ^{m}_{0},
\end{equation}
where $\langle \cdot, \cdot\rangle$ in the integrand is the $Sp(2m;
{\mathbb R})$-invariant inner product on $V_{\rho }$.  Strictly
speaking, this duality has no analog for tensor fields with formal
coefficients but we use it to formally extend the notion in order to
define the following $1$-{\it dual} and $2$-{\it dual} spaces of the
space $T(\rho_1) \otimes T(\rho_2)$, as the duality with respect to the
first (or second) factor.

Clearly, if $B: T(\rho _1) \otimes T(\rho _2)\longrightarrow T(\tau)$
is a $\pbf{\Diff_{\omega}(M)}$-invariant differential operator, then
the operators $B^{*1}: T(\tau)\times T(\rho _2) \longrightarrow T(\rho
_1)$ and $B^{*2}: T(\rho _1)\times T(\tau) \longrightarrow T(\rho
_2)$, the 1-dual and 2-dual of $B$ with respect to the inner product
(IP), are also differential and invariant ones.

The isomorphisms between various realizations of $T(\rho)$ spoken
about in Remark~1 are, clearly, 0-th order invariant differential
operators.  Our description of invariant operators is given up to such
isomorphisms. For the classification of binary operators invariant
with respect to the group of general diffeomorphisms see~\cite{G1}, for
preliminary results on $\Diff_{\omega}(M)$-invariant differential
operators see~\cite{G2}. The results of this paper were preprinted
in~\cite{G3}.

{\bf 0-th order operators} are of the form
\[
Z(\chi , \theta ) = \pr(\chi (x) \otimes \theta (x)),
\]
where $\pr : V_{\rho } \otimes V_{\sigma }\longrightarrow
V_{\tau}$ is the projection of the tensor product onto an irreducible
component.

\medskip

{\bf 1-st order operators} are given by the following theorem

\medskip

\noindent
{\bf Theorem.} {\it Any bilinear $1$-st order (with respect to all
arguments) $\Diff_{\omega }(M)$-invariant differential operator $B:
T(\lambda) \otimes T(\mu)\longrightarrow T(\nu)$ is a linear
combination of the following cases {\em P1--P8} (some of which host
several distinct operators being restricted onto tensors with 
irreducible fibers) and the operators obtained from them by
$1$-dualization or $2$-dualization or transposition of the arguments.}

\medskip

\begin{itemize}
\topsep0mm \partopsep0mm \parsep0mm \itemsep0mm \item[P1)] $\lambda =
(\underbrace{1, \dots , 1}_{p\text{-many\; 1's}}, 0, \dots , 0) $;
weights $\mu $ and $\nu $ differ from each other by a unit in $r$
places, $r \equiv p+1\mod 2$.  For $r \le p+1$ there exists a
representation of these operators in the form $Z(d_{+}\omega , \theta
)$ and for $r \le p-1$ there exists a representation of these
operators in the form $Z(d_{-}\omega , \theta )$.  \item[P2)] The Lie
derivative being restricted onto $Sp(2m;{\mathbb R})$-irreducible
subspaces splits into several operators of the form $Z(d_{\pm }\omega,
\theta )$ and an operator
\[
L : \Vect \times T(\rho )\longrightarrow T(\rho )
\]
which cannot be reduced to operators of the form P1).
\end{itemize}

\medskip

\noindent
{\bf Remark 3.} Observe that if $\xi \in {\mathfrak{h}} (M) \subset {\mathfrak{vect}} (M)$ is
a Hamiltonian vector field, then, by identifying ${\mathfrak{h}}(M) $ with $d\Omega
^{0}$, we see that $d_{+}\xi = d_{-}\xi = 0$ and in this case $L$
coincides with the Lie derivative.  Therefore, $L$ determines a
representation of the Lie algebra ${\mathfrak{h}} (M)$ in the space $T(\rho )$.
It is not difficult to show that the invariance of $B$ is equivalent
to its ${\mathfrak{h}} (M)$-invariance:
\[
L(\xi, B(\chi, \theta )) = B(L(\xi, \chi ), \theta ) + B(\chi, L(\xi, \theta ))
\]
for any $\chi \in T(\rho )$, $\theta \in T(\sigma )$, $\xi \in {\mathfrak{h}} (M)$.

\begin{itemize}
\topsep0mm
\partopsep0mm
\parsep0mm
\itemsep0mm
\item[P3)] $S^{k}{\mathfrak{vect}}
 \times S^{l}{\mathfrak{vect}} \longrightarrow S^{k+l-1}{\mathfrak{vect}}$
(clearly, $S^{k}{\mathfrak{vect}} \cong T(k, 0, \dots , 0)$) is the Poisson
bracket (a.k.a.  the {\it symmetric Schouten's concomitant}) on
(polynomial in momenta) functions on $T^{*}M$.
\item[P4)] $\lambda$, $\mu$, $\nu$ are vectors of the form $(2, 1, \dots , 1, 0, \dots
, 0)$ each, with $p$, $q$ and $ r$ non-zero coordinates,
respectively, such that $p+q+r \equiv 0\mod 2$, $|p-q| \le r \le
p+q$, and $p+q+r \le 2m+2$.
\end{itemize}

If all inequalities are strict, then there exist {\it four distinct
operators} defined on the spaces of such fields, otherwise there exist
only {\it two} distinct operators.  For $p+q+r \le 2m$ two of these
four or two operators are obtained as restrictions of the {\it
Nijenhuis bracket}, or its conjugates, onto the subspaces
\[
T(2, 1, \dots , 1, 0, \dots ,0) \subset\Omega^{p}
\otimes_{C^{\infty}(M)}{\mathfrak{vect}}.
\]

\noindent
{\bf Remark 4.} The remaining two operators (i.e., the ones which are
not the restrictions of the Nijenhuis bracket) are new.  I do not know
anything about them except that they exist and the same applies to the
following two cases P5) and P6).

\begin{itemize}
\topsep0mm
\partopsep0mm
\parsep0mm
\itemsep0mm
\item[P5)] $\lambda$, $\mu $ are of the same form as for P4), $\nu = (3, 1, 1,
\dots , 1, 0, \dots , 0)$.  There exists one operator for $|p-q|+1 \le
r \le p+q-1$, $p+q+r \equiv 1\mod 2$, $p+q+r
\le 2m+1.$
\item[P6)] $\lambda $, $\mu $ are the same as in 4), $\nu = (2, 2, 1, \dots ,
1, 0, \dots , 0)$ with $r$ non-zero entries.  The operator exists
under the same conditions on $p$, $q$, $r$ as for P5).
\item[P7)] $\nu = (1, \dots , 1, 0, \dots , 0)$; whereas $\lambda$, $\mu $ and
conditions on $p$, $q$, $r$ are the same as in 5).  In this case there
exists a unique operator which is not reducible to operators of the
form $d_{\pm }Z.$ It is a restriction of the Nijenhuis bracket.
\item[P8)] $\lambda = (2, 0, \dots , 0)$; whereas $ \mu $ and $\nu $ differ
from each other by a unit at one place.  There exists a unique such
operator.  Further on I will give arguments which enable one to express
it, in principle, explicitly.
\end{itemize}

\noindent
{\bf 2nd order operators:}

I could not {\it classify} such
operators so far.  However, I was lucky to find one new invariant
operator, denoted in the literature $Gz$:
\[
Gz : T(2, 0, \dots , 0) \times T(2, 0, \dots , 0)\longrightarrow T(2,
0, \dots , 0).
\]

For $m = 1$ I got the explicit expression for the operator $Gz$ in
1976.  Let me reproduce it.  In coordinates $x, y$ we have $ \omega =
dx \wedge dy$.  Then
\[
\ba{l}
\ds
 Gz : a\cdot dx^{2} + 2b\cdot dxdy + c\cdot dy^{2}, a^\prime \cdot
dx^{2} + 2b^\prime \cdot dxdy + c^\prime \cdot dy^{2}
\vspace{3mm}\\
\ds \phantom{Gz:} \mapsto
\frac{\partial ^{2}g}{ \partial x^{2}}dx^{2} + 2\frac{\partial ^{2}g}{
\partial x\partial y}dxdy + \frac{\partial ^{2}g}{\partial y^{2}}
dy^{2} + (\{c, a^\prime \}-\{a, c^\prime \})dxdy
\vspace{3mm}\\
\ds \phantom{Gz:\mapsto} + \left (\frac{\partial ^{2}a}{\partial y^{2}} - 2\frac{\partial
^{2}b}{\partial x\partial y} + \frac{\partial ^{2}c}{\partial
x^{2}} \right)
\left(a^\prime dx^{2} + 2bdxdy + c^\prime  dy^{2}\right) + \left(\{a, b^\prime
\}-\{b, a^\prime \}\right)dx^{2}
\vspace{3mm}\\
\ds \phantom{Gz: \mapsto} + \left (\frac{\partial ^{2}a}{ \partial y^{2}} - 2\frac{\partial
^{2}b}{\partial x\partial y} + \frac{\partial ^{2}c}{\partial
x^{2}} \right)\left(adx^{2} + 2bdxdy + cdy^{2}\right)
\left(\{b, c^\prime \}-\{c, b^\prime \}\right)dy^{2},
\ea
\]
where $g = ac^\prime - 2bb^\prime + ca^\prime $ and $\{\cdot, \cdot\}$
is the Poisson bracket. An explicit form of $Gz$ for $m>1$ is to be found.

\medskip

\noindent
{\bf Sketch of the proof of the Theorem:}

Set $y_{i} =
x_{2m+1-i}$ $(1\le i\le m)$, $\partial _{i} = \frac{\partial}{\partial
x_{i}}$, $\delta _{i} = \frac{\partial }{\partial y_{i}}$ .  Denote the
elements of the Lie algebra ${\mathfrak{sp}} (2m; {\mathbb R} ) \subset
{\mathfrak{h}} (M)$ by
\[
e^{ii} = y_{i}\partial _{i}, \qquad e_{ii} = x_{i}\delta _{i}, \qquad
\qquad e^{i}_{j} = x_{j}\partial _{i} + y_{i}\delta _{j}.
\]
Then, clearly,
\[
e^{ij} = e^{ji} = y_{i}\partial _{j} + y_{j}\partial _{i}\qquad
\mbox{and}\qquad e_{ij} = e_{ji} = x_{i}\delta _{j} + x_{j}\delta _{i}\qquad
\text{for} \quad i \neq j.
\]

Let $I(\rho )$ be the space of differential operators from $T(\rho )$ into
$C^{\infty }(M)$ with constant coefficients, i.e.,
\[
I(\rho ) = \left\{\sum P_{i}(\partial , \delta )u_{i}\mid u_{i} \in
V^{*}_{\rho} \cong V_{\rho }\right\}.
\]

The grading in $I(\rho )$ is induced by that in the space of polynomials
$P_{i}$'s,  i.e., $I(\rho )_{0} \cong V_{\rho }$. Define the pairing $I(\rho )
\times T(\rho )\longrightarrow {\mathbb R} $ by the formula
\[
\langle Pu, x\rangle = P(\langle u, \chi (x)\rangle)|_{x=0} .
\]

On $I(\rho )$, define the ${\mathfrak{h}} (M)$-action, dual to the action on
$T(\rho )$, via $L.$ Now, to describe the invariant operators it
suffices to find all the ${\mathfrak{h}} (M)$-morphisms $I(\tau)\longrightarrow
I(\rho ) \otimes _{{\mathbb R}} I(\sigma )$.  It turns out that such a
morphism is completely defined by the image of the highest vector $v
\in V_{\tau} = I(\tau)_{0}$.  Here we have fixed a Borel subalgebra
$\Big\{ \sum\limits_{i\le j}a_{ij}x_{j}\partial _{i}\Big\} \bigcap {\mathfrak{sp}}
(2m;{\mathbb R} )$ so that $w \in I(\rho ) \otimes _{{\mathbb R} }I(\sigma )$ can be
the image of a highest weight singular vector if and only if
\[
\ba{l}
\ds e^{i}_{i+1}w = 0\quad  \text{for}\quad 1 \le i \le m-1\qquad \text{and}
\vspace{2mm}\\
\ds e^{m, m}w = 0 \qquad
(\text{conditions on $w$ to be a highest vector})
\ea
\]
and
\[
\left(x^{2}_{1}\delta _{1}\right)w = 0 \qquad (\text{conditions of {\it
singularity} of the vector})
\]

The degree of $w \in I(\rho ) \otimes _{{\mathbb R} }I(\sigma )$ is equal to
the order of the corresponding differential operator.  The general
form of a vector of degree 1 is
\[
w = \sum^{}_{i\le m} \partial ^{^\prime }_{i}z^{0}_{i} +
\delta ^{^\prime }_{i}t^{0}_{i} + \partial
^{^{\prime\prime}}_{i}z^{1}_{i} + \delta ^{^{\prime\prime}}_{i}t^{1}_{i},
\]
 where $z^{j}_{i}, t^{j}_{i} \in V_{\rho }\otimes V_{\sigma }$,
$ \partial ^\prime (u\otimes v) = \partial u\otimes v$,
$\partial ^{\prime\prime}(u\otimes v) = u\otimes \partial v$. If $w$
 is a highest vector, then all vectors $z, t$ are expressed in terms
 of $z^{0}_{1}$,  $z^{1}_{1}$ which should satisfy
\[
e^{i}_{i+1}z^{j}_{1} = 0\qquad  \text{for}\quad  2 \le i \le m-1,
\qquad \left(e^{1}_{2}\right)^{2}z^{j}_{1} = 0, \qquad e^{m, m}z^{j}_{1} = 0.
\]

The condition $\left(x^{2}_{1}\delta _{1}\right)w = 0$ is  equivalent  to
the  equation
\[
e^{^\prime }_{1, 1}z^{0}_{1} + e^{^{\prime\prime}}_{1, 1}z^{1}_{1} = 0,
\]
where (double) prime means that the operator acts only on the first
(second) multiple of the tensor product.

I have succeeded to define all the cases, where the above system
possesses a solution in $V_{\rho }\otimes V_{\sigma }$; though in
certain cases I was not able to find the solution itself.

Here is an example of a successfully solved case (case 8)):
\[
\lambda = (2, 0, \dots  , 0), \qquad
 \nu = (\mu _{1}, \dots  , \mu _{k-1}, \mu _{k}+1, \mu _{k+1}, \dots);
\]
the case $\nu _{k} = \mu _{k}-1$ is dual to
this one. Let $u_{0} \in V_{\rho }$ be a highest vector, then
\[
u_{0}\otimes v - \frac{1}{ 2}\sum^{}_{2\le i\le k}
e^{i}_{1}u_{0} \otimes e^{1}_{i}v \in V_{\rho }\otimes V_{\sigma }
\]
is a highest vector of weight $(\nu _{1}+1, \nu _{2}, \dots , \nu
_{m})$.  We conclude that
\[
\ba{l}
\ds  z^{0}_{1} = u_{0}\oplus e_{11}v - \sum_{2\le i\le
k}e^{i}_{1}u_{0}\otimes e_{11}e^{1}_{i}v - \frac{1}{ 2}
\sum_{2\le i<j\le m}e^{i}_{1}e^{j}_{1}u_{0} \otimes \left(e_{ij} +
e_{1j}e^{i}_{1} + e_{1i}e^{1}_{j}\right)v
\vspace{3mm}\\
\ds \phantom{z^{0}_{1} =} - \frac{1}{ 2} \sum^{}_{2\le i\le k}(e^{i}_{1})^{2}u_{0} \otimes
 \left(e_{ii} + e_{1i}e^{1}_{i}\right)v
\vspace{3mm}\\
\ds \phantom{z^{0}_{1} =} +
\frac{1}{ 2} \sum_{2\le i\neq j\le m}
e^{i}_{1}e_{1j}u_{0}\otimes \left(e^{j}_{i} + e^{j}_{1}e^{1}_{i}\right)v +
\left(\nu _{i} + e^{i}_{1}e^{1}_{i}\right)v
\vspace{3mm}\\
\ds \phantom{z^{0}_{1} =} +
 \frac{1}{ 4}\sum_{2\le i\le m}e_{ii}\left(e^{i}_{1}\right)^{2}u_{0}\otimes +
 \frac{\nu _{1}-1}{ 2}\sum_{2\le i\le k}e_{11}e^{i}_{1}u_{0}\otimes
 e^{1}_{i}v,
\ea
\]
and
\[
z^{1}_{1} = -e_{11}u_{0}\otimes v + \sum_{2\le i\le
k}e_{11}e^{i}_{1}u_{0}\otimes e^{1}_{i}v.
\]

\noindent
{\bf Conjectures:}
\begin{itemize}
\topsep0mm
\partopsep0mm
\parsep0mm
\itemsep0mm
\item[1)] {\it The operator $Gz$ is a particular case of a
more general operator}:
\[\hspace*{-8mm}
Gz_{r, s}: T(2, \underbrace{1, \dots , 1}_{r\text{-many}}, 0, \dots ,
0) \times T(2, \underbrace{1, \dots , 1}_{s\text{-many}}, 0, \dots ,
0)\longrightarrow T(2, \underbrace{1, \dots , 1}_{(r+s)\text{-many}}, 0,
\dots , 0).
\]
\item[2)] {\it Operators of order $>2$ are compositions of operators of
orders $\leq 2$. There are no operators of order $>5$}.
\end{itemize}

\medskip

\noindent
{\bf Acknowledgements.} I am thankful to A~Kirillov for raising the problem and
D~Leites for general help. Financial support of Swedish Institute is
gratefully acknowledged.

\label{grozman-lastpage}


\newpage

\begin{thebibliography}{99}
\small
\topsep0mm
\partopsep0mm
\parsep0mm
\itemsep0mm

\bibitem{G1}
Grozman P, Classification of Bilinear Invariant Differential
Operators on Tensor Fields, {\it Funct.  Anal.  Appl.}, 1980, V.14, N~2,
58--59.

\bibitem{G2}
Grozman P, The Local Invariant Bilinear Operators on Tensor
Fields on the Plane, {\it Vestnik MGU}, 1980, N~6, 3--6.

\bibitem{G3}
Grozman P, On Bilinear Invariant Differential Operators Acting on
Tensor Fields on a Symplectic Manifold,  in  Seminar
on Supermanifolds, Editor D~Leites, {\it Reports of Dept.  of Math.  Univ. of
Stockholm}, 1988, N~8, 3.

\bibitem{K1}
Kirillov A A, Natural Differential Operations on Tensor Fields,
 Preprint N~56, Akad.  Nauk SSSR,  Inst.  Prikl.  Mat., 1979, 28~p. (in Russian).

\bibitem{K2}
Kirillov A A, Invariant Operators over Geometric Quantities,
in Current Problems in Mathematics, Akad.  Nauk SSSR, Vsesoyuz.
Inst.  Nauchn.  i Tekhn. Informatsii, Moscow, 1980, V.16, 3--29, 228 (in Russian).

\bibitem{L1}
Leites D, Lie Superalgebras,  {\it JOSMAR}, 1984, V.30, N~6, 2481--2513.

\bibitem{LKW}
Leites D, Kochetkov Yu and Vaintrob A, New Invariant Differential
Operators on Supermanifolds and Pseudo-(Co)Homology, in
 General Topology and its Applications,
Editors S~Andima et. al., Marcel Decker, NY, {\it LN in Pure and
Applied Math.}, 1991, V.134, 217--238.



\bibitem{R1}
Rudakov A N, Irreducible Representations of Infinite Dimensional
Lie Algebras of Cartan Type, {\it Math.  USSR Izvestiya}, 1974, V.38,
N~4, 835--866.

\bibitem{R2}
Rudakov A N, Irreducible Representations of Infinite Dimensional
Lie Algebras of Types $S$ and $H$, {\it Math.  USSR Izvestiya}, 1975,
V.39,  N~3, 496--511.

\end{thebibliography}
\end{document}